\documentclass[a4paper,10pt]{article}
\usepackage{amsmath}
\usepackage{amssymb}
\usepackage{graphics}
\usepackage{rotating}
\usepackage{cite}
\usepackage{color}

\DeclareMathOperator{\arctanh}{arctanh}

\title{Cullen-regular quaternionic functions in a Fueter operator framework}
\author{Daniel Alay\'{o}n-Solarz}

\begin{document}

\maketitle

\begin{abstract}
We show characterizations of the class of Cullen-regular functions in the sense of Gentili-Struppa for any domain $\Omega$ in terms of the Fueter operator. We then state a Integral Theorem and discuss how it can be used to define a more general version of Cullen-regularity, that does not require the function to be of class $C^{1}$.
\end{abstract}

\section{Preliminaries}
Let \textbf{H} be the algebra of the quaternions and let $p$ be a quaternion, then $p$ can be written as
\begin{displaymath}
p := t + xi + yj + zk.
\end{displaymath}
Define $\iota$ as
\begin{displaymath}
\iota := \frac{xi+yj+zk}{\sqrt{x^{2}+y^{2}+z^{2}}}.
\end{displaymath}
then we write $p$ as $t + r \iota$, with $r$ the norm of the imaginary part of $p$. \\
Recall the left-Fueter operator is given by: 
\begin{displaymath}
 D_{l}:=\frac{\partial  }{\partial t}+ i\frac{\partial }{\partial x} + j\frac{\partial }{\partial y} + k\frac{\partial }{\partial z}. 
\end{displaymath}
Note that $\iota$ can be parametrized by spherical coordinates
\begin{displaymath}
\iota = (\cos\alpha \sin\beta, \sin\alpha \sin\beta, \cos \beta).
\end{displaymath}
The coordinate system based in the variables $(t,r,\alpha,\beta)$ is especially well suited to study the interplay between the Fueter operator and the Cullen-regular quaternionic functions. The Fueter operator in this coordinate system has the form:
\begin{displaymath}
D_{l} = \frac{\partial }{\partial t} + \iota \frac{\partial }{\partial r}  - \frac{1}{r} \frac{\partial}{\partial_{l} \iota},
\end{displaymath}
for all $p$ that does not lie in the complex subplane of the quaternions given by $t + zk$. Where
\begin{equation*}
\frac{\partial}{\partial_{l} \iota} := ({\iota}_{\alpha})^{-1}\frac{\partial}{\partial \alpha} +  ({\iota}_{\beta})^{-1}\frac{\partial}{\partial \beta}.
\end{equation*}
and $\iota_{\alpha}$ and $\iota_{\alpha}$ represent the derivatives of $\iota$ with respect to $\alpha$ and $\beta$ respectively.
Quaternionic functions of one quaternionic variable that are of class $C^{1}$ and null-solutions to the \textit{left-Cullen operator}, defined as:
\begin{displaymath}
 \frac{\partial }{\partial t} + \iota \frac{\partial }{\partial r} 
\end{displaymath}
will be called \textit{left-Cullen regular}.  We denote the real numbers in the quaternions as \textbf{R}. 
\section{A characterization of Cullen-regular functions}

\newtheorem{theorem}{Theorem}
\begin{theorem}
Let $f$ be a $C^{1}$ quaternionic function of one quaternionic variable  $p$ with domain $\Omega$,  the following statements are equivalent:
\begin{enumerate}
\item { $f$ is left-Cullen-regular in $\Omega$.}
\item $f$ satisfies
\begin{displaymath}
D_{l} \iota f + \iota D_{l}  f = \frac{-2f}{r}
\end{displaymath}
in $\Omega \backslash \textbf{R}$
\item
There exists quaternionic functions $u$ and $v$ of one quaternionic variable $p$ such that $f= u +\iota v$  in $\Omega$ and \begin{displaymath}
D_{l}f = \frac{-2v}{r}
\end{displaymath}
and
\begin{displaymath}
D_{l}(\iota f) = \frac{-2u}{r}
\end{displaymath}
in $\Omega \backslash \textbf{R}$
\item
With the same $u$ and $v$ that above:
\begin{displaymath}
D_{l}(\frac{f}{r^2}) = \frac{-2 \iota u }{r^{3}}
\end{displaymath}
and
\begin{displaymath}
  D_{l}(\frac{\iota f}{r^2}) = \frac{-2 \iota v }{r^{3}}
\end{displaymath}
\end{enumerate}
\end{theorem}

The proof is based on calculations using the given expression for the Fueter operator.  Is important the following fundamental property of the angular part of the Fueter operator.
\newtheorem{lemma}{Lemma}
\begin{lemma}
For every quaternionic function $f$ differentiable in the variables $\alpha$ and $\beta$ it holds:
\begin{displaymath}
\frac{\partial }{\partial\iota}(\iota f)+ \iota \frac{\partial }{\partial\iota}(f) = 2f
\end{displaymath}
\end{lemma}
And the following choices for $u$ and $v$.
\begin{displaymath}
u := \frac{1}{2}\frac{\partial  }{\partial\iota}(\iota f)
\end{displaymath}
and
\begin{displaymath}
v := \frac{1}{2}\frac{\partial }{\partial\iota}(f)
\end{displaymath}
together with the observation that a function $f$ is left-Cullen-regular in the domain $\Omega$ if and only if $\iota f$ is left-Cullen-regular on their common domain $\Omega \backslash \textbf{R}$.

\section{Special cases}
1) The domain $\Omega$ of the Cullen-regular function $f$ is a ball of radius $R$ centered on the origin or a spherical shell with radius $R_{1}, R_{2}$. Then $f$ admits a power series and Laurent series representation with quaternionic coefficients on the right, see \cite{GS2}.  For functions $f$ defined as such power series or Laurent series, as a consequence of the Fueter's Theorem one has:
\begin{displaymath}
D_{l} \Delta f = 0.
\end{displaymath}
Note that this representation fails if we don't assume the function is defined on some open subset of the real numbers. One example is the Cullen-regular function $f := \iota$.\\
2) The functions $u$ and $v$ defined above are real-valued functions (such as the function $f:= \iota$). Then $u$ and $v$ must satisfy the following equations, see \cite{Ala1}:

\begin{eqnarray}
\frac{\partial v}{\partial \alpha}(\sin\beta)^{-1}+\frac{\partial u}{\partial \beta} = 0,
\\
\frac{\partial u}{\partial \alpha}(\sin\beta)^{-1}-\frac{\partial v}{\partial \beta}= 0.
\end{eqnarray}
If the Cullen-regular function $f$ is defined over the whole imaginary $2$-sphere and the  $u$ and $v$ functions are real-valued then they must be constants as functions of $\alpha$ and $\beta$. If the function is defined in almost all of the 2-sphere and $u$ and $v$ are real then there might exists non-trivial solutions. Consider the following 3 Cullen-regular functions of the variable $p$:
\begin{displaymath}
\arctan{\frac{x}{y}} + \iota \arctanh {\frac{z}{r}}
\end{displaymath}
\begin{displaymath}
\arctan{\frac{y}{z}} + \iota \arctanh {\frac{x}{r}}
\end{displaymath}
\begin{displaymath}
\arctan{\frac{z}{x}} + \iota \arctanh {\frac{y}{r}}
\end{displaymath}
3) In the spirit of the previous case, if the functions $u$ and $v$ of the Cullen-regular function $f$ satisfy equations (1) and (2)
as 4-vectors. We refer to this latter case saying that the Cullen-regular function $f$  is \textit{Hyperholomorphic}. Note that in particular constants functions, real or not, are hyperholomorphic. Because multiplying left-hyperholomorphic functions by quaternionic constants by the right is again a hyperholomorphic function and because the power function $p \to p^n$ is hyperholomorphic for all integers $n$ we see that functions constructed as power series or Laurent series of quaternions with quaternionic coefficients at the right are hyperholomorphic. On the other side, if we add the special assumption that the hyperholomorhic $f$ is 
of class $C^{3}$ it can be asserted that the hyperholomorphic function $f$, without any assumption on $\Omega$ satisfy the Fueter's Theorem 
\begin{displaymath}
D_{l} \Delta f = 0,
\end{displaymath}
even if the function $f$ does not satisfy Fueter's original hypothesis or it does not express as a power series or Laurent series with quaternionic coefficients by the right. By imposing a different condition, for example that the hyperholomorphic function $f$ is such that $u$ and $v$ are real valued, then it can be asserted that for any hyperholomorphic function $g$ the function $fg$ is hyperholomorphic in their common domain. Which is immediate if $f$ and $g$ can be expressed as quaternionic power series or Laurent series with coefficients by the right and $f \iota  = \iota f$. \\
With appropiate changes analog results can be expressed for right-Cullen regular functions.

\section{A Integral Theorem for Cullen-regular functions}
The Fueter operator approach to Cullen-regular functions allows to state a integral theorem, first we have the following consequence of Gauss Theorem:
\begin{lemma}
Let $(n_{0},n_{1},n_{2},n_{3})$ be the outward unit normal of  $K$,  any smooth, simple closed hypersurface in the quaternionic space and $K^{*}$ its interior, then
\begin{displaymath}
\int_{K}(f_{0}n_{0}+f_{1}n_{1}+f_{2}n_{2}+f_{3}n_{3}) dS_{K} = \int_{K^{*}}(\frac{\partial f_{0}}{\partial t}+\frac{\partial f_{1}}{\partial x}+\frac{\partial f_{2}}{\partial y}+\frac{\partial f_{3}}{\partial z}) dV
\end{displaymath}
where $f_{i}$ are differentiable quaternionic functions.
\end{lemma}
from this we obtain
\begin{theorem}
Let $(n_{0},n_{1},n_{2},n_{3})$ be the outward unit normal of  $K$,  any smooth, simple closed hypersurface in the quaternionic space , disjoint from the real axis and $K^{*}$ its interior, if $f$ is left-Cullen-regular
\begin{displaymath}
\int_{K}n(p) f(p) dS_{K} = \int_{K^{*}}(-\frac{2v}{r})dV
\end{displaymath}
where $dS_{K}$ is the element of surface area on $K$.
\end{theorem}
Suppose $f$ is not $C^{1}$, then $f$ is not Cullen-regular in the classic sense. Now suppose $f$ is continous and admits derivatives as a function of $\alpha$ and $\beta$. This function might not be Cullen-regular in the classic sense. However in this case one sees the integrals of the left and right hand side of the Integral Theorem exists. This motivate us to define a generalized left-Cullen-regular function as a continuous function that is a $C^{1}$ as function of $\alpha$ and $\beta$ defined on a domain $\Omega$  such that the integral theorem \textit{and} its $\iota f$ counterpart hold for all smooth, simple closed hypersurfaces $K$ disjoint from the real axis. Clearly classic Cullen-regular functions are also Cullen-regular in this generalized sense.

%%%%%%%%%%%%%%%%%%%%%%%%%%%%%%%%%%%%%%%%%%%%%%%%%%%%%%%%%%%%%
%                                SECTION V
%%%%%%%%%%%%%%%%%%%%%%%%%%%%%%%%%%%%%%%%%%%%%%%%%%%%%%%%%%%%%%%%%%%%%%%%%%%%%%%

%%%%%%%%%%%%%%%%%%%%%%%%%%%%%%%%%%%%%%%%%%%%%%%%%%%%%%%%%%%%%%%%%%%%%%%%%%%%%%%
%                                SECTION VI
%%%%%%%%%%%%%%%%%%%%%%%%%%%%%%%%%%%%%%%%%%%%%%%%%%%%%%%%%%%%%%%%%%%%%%%%%%%%%%%

%%%%%%%%%%%%%%%%%%%%%%%%%%%%%%%%%%%%%%%%%%%%%%%%%%%%%%%%%%%%%%%%%%%%%%%%%%%%%%%
%                                SECTION VII
%%%%%%%%%%%%%%%%%%%%%%%%%%%%%%%%%%%%%%%%%%%%%%%%%%%%%%%%%%%%%%%%%%%%%%%%%%%%%%%

%%%%%%%%%%%%%%%%%%%%%%%%%%%%%%%%%%%%%%%%%%%%%%%%%%%%%%%%%%%%%%%%%%%%%%%%%%%%%%%
%                                REFERENCES
%%%%%%%%%%%%%%%%%%%%%%%%%%%%%%%%%%%%%%%%%%%%%%%%%%%%%%%%%%%%%%%%%%%%%%%%%%%%%%%

\end{document}